\documentclass[a4paper,12pt,final]{amsart}
\usepackage{times,a4wide,mathrsfs}

\newcommand{\C}{\mathbb{C}}

\newtheorem{theorem}{Theorem}
\newtheorem{lemma}{Lemma}

\newtheorem{proposition}{Proposition}
\newtheorem{remark}{Remark}

\begin{document}

\title{A sharp  bound on the Lebesgue constant for Leja points in the unit disk}

\author[M. Ouna\"{\i}es]{Myriam Ouna\"{\i}es}
\address{Institut de Recherche Math\'ematique Avanc\'ee,
Universit\'e de Strasbourg, 7 Rue Ren\'e Des\-car\-tes,
67084 Strasbourg CEDEX, France}
\email{ounaies@math.u-strasbg.fr}

\subjclass[2000]{}

\date{\today}

\begin{abstract}
We give a sharp bound for the Lebesgue constant associated to Leja sequences 
 in the complex unit disk, confirming a conjecture made by Calvi and Phung \cite{Ca-Ph}.
\end{abstract}
\maketitle

\section{Introduction and statement of the main results}

Let $K$ be a compact set 
in the complex plane. Leja sequences  are defined by   arbitrarily fixing a first point $e_0\in K$  then by recursively selecting  $e_k$, $k=1,2,\cdots$  such that 
\begin{equation}\label{Leja}
\prod_{j=0}^{k-1}\vert e_k-e_j\vert=\max_{z\in K} \prod_{j=0}^{k-1}\vert z-e_j\vert.
\end{equation}
They were first studied by Erdei \cite[page 78]{Er}, then by Leja \cite{Le} who showed that the sequence $\left(\prod_{j=0}^{k-1}\vert e_k-e_j\vert\right)^{\frac{1}{k}}$ converges to the tranfinite diameter of $K$.

Consider ${\mathcal C}({K})$ the space of continous functions on $K$, endowed with the uniform norm. For any $f\in {\mathcal C}({K})$, the unique polynomial in the space $\Pi_{k-1}$ of polynomials of degree at most $k-1$ which coincides with $f$ on $E_k=(e_0,e_1,\cdots,e_k)$ is the  Lagrange interpolation polynomial defined by
\[L_{E_k}(f)(z)=\sum_{j=0}^{k-1}f(e_j)l_{j,E_{k}}(z)\]
where 
\[ l_{j,E_{k}}(z)=\prod_{i=0,i\not=j}^{k-1}\frac{z-e_i}{e_j-e_i}\] are the fundamental Lagrange interpolation polynomials.

The norm of $L_{E_{k}}$ as a a continous linear operator  from ${\mathcal C}({K})$ into $\Pi_{k-1}$ is the so-called Lebesgue constant 
\begin{equation}\label{Lebesgue}
\Lambda_{E_k}:=\sup_{\Vert f\Vert\le 1}\Vert L_{E_k}(f)\Vert=\sup_{z\in K}\vert\lambda_{E_k}(z)\vert,
\end{equation}
where  
\[\lambda_{E_k}(z):=\sum_{j=0}^{k-1}\vert  l_{j,E_{k}}(z)\vert.\]

The Lebesgue constant  plays a crucial role in polynomial interpolation.
The inequality
\[\Vert L_{E_k}(f)-f\Vert_K\le (\Lambda_{E_k}+1)\inf_{p\in \Pi_{k-1}}\Vert f-p\Vert,\]
shows that it measures how close the interpolant is to the best polynomial approximant of a function.
 It also measures the stability of the Lagrange interpolation. We refer the reader to \cite{Ve} for more details and many interesting properties on the Lebesgue constant.

\vskip0.2cm
 In this paper, we are interested in finding an optimal bound for the Lebesgue constant associated with  Leja points in  the case where $K$ is the complex unit disk $\mathcal U=\{z\in \C\ :\ \vert z\vert\le 1\}$. 

We consider  Leja sequences $(e_k)_{k\ge 0}$ initiated at $e_0\in \partial {\mathcal U}=\{\vert z\vert=1\}$. There is no loss of generality in assuming that $e_0=1$ since any Leja sequence is the product by $e_0$ of a Leja sequence initiated at $1$.  We speak of Leja sequences rather than of a Leja sequence because at each step there may be more than one point $e_k$ satisfying ($\ref{Leja}$).
Note that by the maximum principle, $\vert e_k \vert=1,\ \ k=1,2,\cdots$. 

Any finite sequence $E_k=(e_0,\cdots,e_k)$ where $e_0,\cdots,e_k$ are defined by $(\ref{Leja})$ is called a $k$-Leja section.

 Leja sequences of the disk were explicitely described by Bialaz-Ciez and Calvi  \cite{Bi-Ca}. They showed   that for any  Leja sequence of the disk, initiated at $e_0=1$, the underlying set of the $2^n$-th section consists of the $2^n$-th roots of the unity. The $2^{n+1}$-th section is 
\[(E_{2n},\rho F_{2n})\]
where $\rho$ is any $2^n$-th root of $-1$ and $F_{2^n}=(f_0,\cdots,f_{2^n-1})$ is a $2^n$-Leja section  with $f_0=1$,
with the notation
\[(E_k,F_j):=(e_0,\cdots,e_{k-1},f_0,\cdots,f_{j-1}).\]

\noindent A natural Leja sequence is then inductively constructed as follows :
\[E_1=(1),\  \ (E_{2^n}, e^{\frac{i\pi}{2^n}}E_{2^n}).\]

\noindent It was shown in \cite{Bi-Ca} that this particular sequence is defined by 
\begin{equation}\label{Def}
e_k=\exp(i\pi \sum_{j=0}^s  2^{-p_j})
\end{equation}
where $k\ge 1$ is expanded in the binary form
\begin{equation}\label{exp} 
k=2^{p_0}+\cdots+2^{p_{s}}, \ \ 0\le p_0<p_1<\cdots<p_s.
\end{equation} 

\noindent Let us now recall the previous results about the bounds of the Lebesgue constant for Leja points in the disk.

\vskip 0.2cm
Calvi and Phung \cite{Ca-Ph}
  showed, for any $k$-Leja section $E_k=(e_0,\cdots,e_k)$,  that
\[
\Lambda_{E_k}= O(k\ln k).
\]  
They also showed that 
\[ \Lambda_{E_{2^n-1}}=2^n-1,\ \ n\ge 1\] so it cannot grow slower than $k$. 
They conjectured the following :
\begin{equation}\label{conjecture}
\Lambda_{E_k}\le k
\end{equation}

Chkifa \cite{Ch1} gave the sharper bound : 
\[
\Lambda_{E_k}\le 2k.
\]

In the same direction, Irigoyen \cite{Iri} showed that there is a uniform bound for the fundamental Lagrange interpolation polynomials, namely :
\[
\sup_{k> j\ge 0}\left(\sup_{z\in {\mathcal U}} \vert l_{j,E_k}(z)\vert\right)\le \pi\exp(3\pi).
\]

In a recent work Chkifa \cite{Ch2}  introduced the so-called "quadratic" Lebesgue function 
\[ \lambda_{{E_k},2}(z)=\left(\sum_{j=0}^{k-1}\vert l_{j,E_k}(z)\vert^2\right)^{1/2}\]
and the "quadratic" Lebesgue constant 
\[\Lambda_{{E_k},2}:=\max_{z \in {\mathcal U}} \lambda_{{E_k},2}(z),\] which turned out to be  quite efficient tools. Their advantage is that they strongly exploit the binary structure  of the Leja sequences. 
It was proved in \cite{Ch2} that
for any $k$ expanded in the form $(\ref{exp})$, the following inequalities hold :
\begin{equation}\label{Chkifa}
\sqrt{2^{s+1}-1}= \lambda_{{E_k},2}(e_k)\le \Lambda_{{E_k},2}\le \sqrt{3 (2^{s+1}-1)}
\end{equation}
and by the Cauchy-Schwarz inequality, 

\[
\Lambda_{E_k}\le \sqrt{3k (2^{s+1}-1)},
\]

\vskip 0.2cm

 In the present paper, we will follow the same approach as in \cite{Ch2}. Our estimate for the "quadratic" Lebesgue constant is given by next proposition.

\begin{proposition}\label{main}
For any $k\ge 1$ and any  $k$-Leja section $E_k=(e_0,\cdots,e_k)$ in the unit disk, we have :
\[\Lambda_{{E_k},2}\le \sqrt{2^{-p_0}k},\]
where $k=2^{p_0}+\cdots+2^{p_{s}}$ is expanded in the form ($\ref{exp}$).

\end{proposition}
\begin{remark}
This a sharp bound : in the case where $k=2^n-1$,  the proposition combined with the first inequality in 
($\ref{Chkifa}$) shows that we actually have
\begin{equation}
\Lambda_{E_{2^n-1},2}=\sqrt{2^n-1}
\end{equation}

\end{remark}
Now, by Proposition $\ref{main}$ and a straightforward application of the  Cauchy-Schwarz inequality, we settle conjecture $(\ref{conjecture})$ :
\begin{theorem}\label{theorem}
For any $k\ge 1$ and any $k$-Leja section $E_k=(e_0,\cdots,e_k)$ in the unit disk, the following inequality holds :
\[\Lambda_{E_k}\le 2^{-\frac{p_0}{2}}k,\]
where $k=2^{p_0}+\cdots+2^{p_{s}}$ is expanded in the form ($\ref{exp}$).

\end{theorem}

\section{proof of proposition\ref{main}}

It was proved in  \cite[Lemma 2.4]{Ch1} that for any $k\ge 0$ and any two $k$-Leja sections $E_k$ and $F_k$ of the unit disk, there exists $\vert \rho\vert=1$ such that $F_k=\rho E_k$ (in the set sense). This implies that the Lebesgue constant of a $k$-Leja sequence of the disk only depends on $k$.  

We may then assume, without loss of generality, that we are dealing with the particular Leja-sequence, that we will denote by $E=(e_k)_{k\ge 0}$, defined by $(\ref{Def})$. 

We will use the simplified notations :
\[l_{j,k}(z)=\prod_{i=0,i\not=j}^{k-1}\frac{z-e_i}{e_j-e_i},\ \ \ 0\le j\le k-1,\]
\[ \lambda_{k,2}(z)=\left(\sum_{j=0}^{k-1}\vert l_{j,k}(z)\vert^2\right)^2,\]
\[\Lambda_{k,2}(z)=\sup_{\vert z\vert\le 1} \lambda_{k,2}(z)=\sup_{\vert z\vert= 1} \lambda_{k,2}(z).\]
The last equality comes from the subharmonicity of $\lambda_{k,2}$ and the maximum principle.

 Recall that, by definition of Leja sequences,
\begin{equation}\label{1}
\vert l_{k-1,k}(z)\vert=\prod_{j=0}^{k-2}\frac{\vert z-e_i\vert}{\vert e_{k-1}-e_i\vert}\le 1,\ \ \ k\ge 2.
\end{equation}

Exploiting the following relations satisfied by these Leja points :
\begin{equation}\label{Sym}
e_{2j+1}=-e_{2j},\ e_{2j}^2=e_{2j+1}^2=e_j, \ \ \ j\ge 0,
\end{equation} 
it was shown in \cite{Ch2} that :
\[\lambda_{2N,2}(z)=\lambda_{N,2}(z^2),\ \ \ N\ge 0,\ \vert z\vert \le 1\]
and as consequences : 
\[
\Lambda_{2N,2}=\Lambda_{N,2},\ \ \ N\ge 1
\]
\[
\Lambda_{2^n,2}=\Lambda_{1,2}=1,\ \ \ n\ge 0.
\]
To get an accurate  estimate on the Lebesgue constant, we need a  recursive formula for $\Lambda_{{2N+1},2}$. We will use similar techniques as in \cite{Ch2} to show the following :
\begin{lemma}\label{Ineg}
For all  $N\ge 1 $,
\[ \Lambda^2_{2N+1,2}\le \frac{1}{2}\Lambda^2_{N+1,2}+2\Lambda^2_{N,2}+\frac{1}{2}.\]
\end{lemma}

\begin{remark}
We may already deduce by induction from this lemma that, for all $n\ge 1$,
\[\Lambda^2_{2^n-1,2}\le 2^n-1.\]
\end{remark}
\begin{proof}${}$

Let $\vert z\vert\le 1$. 
Using  the relations $(\ref {Sym})$, we find
\[ l_{2N,2N+1}(z)=
\prod_{i=0}^{N-1}\frac{z^2-e_i}{e_N-e_i}=l_{N,N+1}(z^2).\]
For $j=0,1,\cdots,N-1$, we have 
\[
\begin{split}
l_{2j,2N+1}(z)&=
\frac{(z-e_{2N})(z-e_{2j+1})}{(e_{2j}-e_{2N})(e_{2j}-e_{2j+1})}\prod_{i=0,i\not=j}^{N-1}\frac{(z^2-e_i^2)}{(e_j-e_i)}\\
&=\frac{(z-e_{2N})(z+e_{2j})(e_{2j}+e_{2N})}{2e_{2j}(e_j-e_N)}l_{j,N}(z^2).\\
\end{split}
\]
And similarily,
\[
\begin{split}
l_{2j+1,2N+1}(z)
=\frac{(z-e_{2N})(z-e_{2j})(e_{2j}-e_{2N})}{2e_{2j}(e_j-e_N)}l_{j,N}(z^2).\\
\end{split}
\]
Now, we use that :
\[
\vert (z+e_{2j})(e_{2j}+e_{2N})\vert^2+\vert (z-e_{2j})(e_{2j}-e_{2N})\vert^2=2\vert z+e_{2N}\vert^2+2 \vert ze_{2N}+e_j\vert^2\\
\]
and that
\[ \vert((z-e_{2N})(ze_{2N}+e_j)\vert^2\le 2\vert z^2-e_j\vert^2+2\vert e_j-e_N\vert^2.\]
We deduce the following :
 \[
\begin{split}
\vert l_{2j,2N+1}(z)\vert^2+\vert l_{2j+1,2N+1}(z)\vert^2
&\le \frac{\vert z^2-e_N\vert^2+2\vert z^2-e_j\vert^2+2\vert e_j-e_N\vert^2}{2\vert e_j-e_N\vert^2}\vert l_{j,N}(z^2)\vert^2\\
&=\frac{1}{2}\vert l_{j,N+1}(z^2)\vert^2+\vert l_{N,N+1}(z^2)\vert^2\vert l_{j,N}(e_N)\vert^2+\vert l_{j,N}(z^2)\vert^2\\
\end{split}
\]
We are now ready to estimate $\lambda_{2N+1,2}(z)$ :
\[
\begin{split}
\lambda^2_{2N+1,2}(z)
&\le\frac{1}{2}\lambda^2_{N+1,2}(z^2)+\vert l_{N,N+1}(z^2)\vert^2\lambda^2_{N,2}(e_N)+\lambda^2_{N,2}(z^2)+\frac{1}{2}\vert l_{N,N+1}(z^2)\vert^2.\\
\end{split}
\]
Applying inequality $(\ref{1})$, the proof of the lemma is completed.

\end{proof}

Let us now define the sequence $(U_k)_{k\ge 1}$ by :
\[
U_1=1, U_{2N}=U_N,\ U_{2N+1}=\frac{1}{2}U_{N+1}+2U_N+\frac{1}{2}, \ \ \ N\ge 1.\]
We easily verify by induction on $n\ge 0$ that for all $m\ge 0$ :
\begin{equation}\label{ind}
U_{2^nm}=U_m;
\end{equation}
\begin{equation}\label{Gen}
U_{2^nm+1}= 2^{-n}U_{m+1}+4(1-2^{-n})U_m+1-2^{-n}.
\end{equation}
Thanks to  Lemma $\ref{Ineg}$, we have  :
\[
\Lambda^2_{k,2}\le U_k,\ \ k\ge 1.
\]
So the proof of Proposition $\ref{main}$ will be done if we prove that, for all $k= 2^{p_0}+2^{p_1}+\cdots+2^{p_s}$ expanded as in $(\ref{exp})$,
\[
U_k \le 2^{-p_0}k, 
\]
or equivalently, that 
\begin{equation}\label{Delta}
\Delta_k := 2^{-p_0}k-U_k \ge 0.
\end{equation}
Clearly, because of $(\ref{ind})$, we have 
\begin{equation}\label{simp}
\Delta_{2^{p_0}+2^{p_1}+\cdots+2^{p_s}}= \Delta_{1+2^{p_1-p_0}+\cdots+2^{p_s-p_0}}.
\end{equation}
When $s=0$, we note that, for all $p_0\ge 0$ : 
\[
\Delta_{2^{p_0}}=\Delta_1=0.
\]

\noindent Let us now deal with the case $s=1$ :
Thanks to  $(\ref{Gen})$, for all $0\le p_0<p_1$,
\[
U_{1+2^{p_1-p_0}}=2^{p_0-p_1}
U_2\\
+4(1-2^{p_0-p_1})U_1+1-2^{p_0-p_1}=4(1-2^{p_0-p_1})+1.
\]
It follows that
\begin{equation}\label{s1}
\Delta_{2^{p_0}+2^{p_1}}=2^{p_0-p_1}(2^{p_1-p_0}-2)^2\ge 0.\\
\end{equation}
\vskip0.2cm

For the general case $s\ge 2$, $(\ref{Delta})$ will be proved recursively on $s$ in the following  lemma :
 \begin{lemma}
For all $s\ge 1$ and for all $0\le p_0<p_1<\cdots<p_s$ :
\begin{equation}\label{lemma}
\begin{split}
\Delta_{2^{p_0}+\cdots+2^{p_s}}&=
\sum_{j=1}^s 2^{j-p_j+p_0}(2^{p_j-p_{j-1}}-2)\left[(3\ 2^{j-1}-1)\Delta_{2^{p_j}+\cdots+2^{p_s}}\right.\\
&\left.+(2^{j-2}(2^{p_j-p_{j-1}}-2^2)+1)(1+2^{p_{j+1}-p_j}+\cdots+2^{p_s-p_j})\right]\ge 0\\
\end{split}
\end{equation}
\end{lemma}
\begin{remark}${}$

$\Delta_k$ measures how large  $2^{-p_0}k$ is compared to $\Lambda^2_{k,2}$  in terms of the binary expansion of $k$.  

Whenever $p_j-p_{j-1}=1$, the factor $2^{p_j-p_{j-1}}-2$ is a vanishing term and does not contribute to the sum.
This lemma states in particular that :

$\Delta_k=0$  if and only if $k=2^{p_0}(2^n-1)$ for some $p_0\ge 0$ and $n\ge 1$. 

For all the other numbers, the inequalities 
\[\Lambda_k\le \sqrt{k}\Lambda_{k,2}\le \sqrt{k} \sqrt{2^{-p_0}k-\Delta_k}\]
are actually more accurate than the ones announced in Proposition \ref{main} and Theorem \ref{theorem}.
\end{remark}

\begin{proof}${}$

When $s=1$, this corresponds to $(\ref{s1})$.
Assume that $s\ge 2$ and that formula and inequality $(\ref{lemma})$ are verified up to  $s-1$. 

\noindent We apply $(\ref{ind})$ and $(\ref{Gen})$ to get  the next two formulas :
\[
\begin{split}
U_{2^{p_0}+\cdots+2^{p_s}}=2^{p_0-p_1}
U_{2+2^{p_2-p_1}+\cdots+2^{p_s-p_1}}+4(1-{2^{p_0-p_1}})U_{2^{p_1}+\cdots+2^{p_s}}+1-2^{p_0-p_1}.\\
\end{split}
\]
\[
U_{2^{p_1}+2^{p_2}+\cdots+2^{p_s}}=
\frac{1}{2}U_{2+2^{p_2-p_1}+
\cdots+2^{p_s-p_1}}
+2U_{2^{p_2}+\cdots+2^{p_s}}+\frac{1}{2}.
\]
And we deduce that
\[
U_{2^{p_0}+\cdots+2^{p_s}}=1+2^{p_0-p_1}\left[(2^{2+p_1-p_0}-2)U_{2^{p_1}+\cdots+2^{p_s}}-4U_{2^{p_2}+\cdots+2^{p_s}}-2\right].
\]

In terms of the $\Delta_k's$, the previous identity may be written as follows :
\begin{equation}\label{R}
\begin{split}
\Delta_{2^{p_0}+\cdots+2^{p_s}}=&
2^{p_0-p_1}(2^{p_1-p_0}-2)\left[4\Delta_{2^{p_1}+\cdots+2^{p_s}}+(2^{p_1-p_0}-2)(1+2^{p_2-p_1}+\cdots+2^{p_s-p_1})\right]\\
&+2^{p_0-p_1} R_1(p_1,\cdots,p_s),\\
\end{split}
\end{equation}
where
\[
\begin{split}
R_1(p_1,\cdots,p_s):=&-2(2^{p_2-p_1}-2)\left[(1+2^{p_3-p_2}+\cdots+2^{p_s-p_2})+2^{1+p_1-p_2}\Delta_{2^{p_2}+\cdots+2^{p_s}}\right]\\
&+6\Delta_{2^{p_1}+\cdots+2^{p_s}}-2^{3+p_1-p_2}\Delta_{2^{p_2}+\cdots+2^{p_s}} \\
\end{split}
\]
After applying the recursive hypothesis $(\ref{lemma})$ for ranks $s-1$ and $s-2$ on the second line, and setting apart the term corresponding to $j=1$ in $\Delta_{2^{p_1}+\cdots+2^{p_s}}$,  we obtain :
\begin{equation}\label{R1}
\begin{split}
R_1(p_1,\cdots, p_s):=&2^{2+p_1-p_2}(2^{p_2-p_1}-2)\left[5\Delta_{2^{p_2}+\cdots+2^{p_s}}+(2^{p_2-p_1}-3)(1+2^{p_3-p_2}+\cdots+2^{p_s-p_2})\right]\\
&+R_2(p_2,\cdots, p_s)\\
\end{split}
\end{equation}
where
\[
\begin{split}
R_2(p_2,\cdots, p_s)=&3\sum_{j=2}^{s-1}2^{1+j-p_{j+1}+p_1}(2^{p_{j+1}-p_j}-2)\left[(3\ 2^{j-1}-1)\Delta_{2^{p_{j+1}}+\cdots+2^{p_s}}\right.\\
&\left.+(2^{j-2}(2^{p_{j+1}-p_j}-2^2)+1)(1+2^{p_{j+2}-p_{j+1}}+\cdots+2^{p_s-p_{j+1}})\right]\\
&-\sum_{j=1}^{s-2}2^{3+j-p_{j+2}+p_1}(2^{p_{j+2}-p_{j+1}}-2)\left[(3\ 2^{j-1}-1)\Delta_{2^{p_{j+2}}+\cdots+2^{p_s}}\right.\\
&\left.+(2^{j-2}(2^{p_{j+2}-p_{j+1}}-2^2)+1)(1+2^{p_{j+3}-p_{j+2}}+\cdots+2^{p_s-p_{j+2}})\right].\\
\end{split}
\]
Writing the last expression under one sum, we obtain :
\begin{equation}\label{R2}
\begin{split}
R_2(p_2,\cdots, p_s)=&
\sum_{j=3}^s 2^{j-p_j+p_1}(2^{p_j-p_{j-1}}-2)\left[(3\ 2^{j-1}-1)\Delta_{2^{p_j}+\cdots+2^{p_s}}\right.\\
&\left.+(2^{j-2}(2^{p_j-p_{j-1}}-2^2)+1)(1+2^{p_{j+1}-p_j}+\cdots+2^{p_s-p_j})\right].\\
\end{split}
\end{equation}

Finally, we put together equations $(\ref{R})$, $(\ref{R1})$ and $(\ref{R2})$ and we have reached formula (\ref{lemma}).

\end{proof}

\end{document}